\documentclass[12pt,leqno]{article}
\title{Some questions on the class group of cyclotomic fields }
\author{Roland Qu\^eme}
\usepackage{amsmath,amsthm,amstext,amsbsy}

\font\mathbb=msbm10

\newcommand{\Q}{\mbox{\mathbb Q}}
\newcommand{\Z}{\mbox{\mathbb Z}}

\newcommand{\modu}{\ \mbox{mod}\ }
\newcommand{\be}{\begin{equation}}
\newcommand{\ee}{\end{equation}}
\newcommand{\bd}{\begin{displaymath}}
\newcommand{\ed}{\end{displaymath}}
\newcommand{\bn}{\begin{enumerate}}
\newcommand{\en}{\end{enumerate}}
\date{2007 march 08}
\begin{document}
\maketitle
\abstract This article deals with a  study of the structure of  the
class group of the cyclotomic field $K=\Q(\zeta_p)$ for $p$ an odd
prime number, starting from Stickelberger relation. The present
state of this work leads me  to set a question  for all the prime
numbers $h\not=p$ which divide the relative class number
$h^-(K)$.
%
\section{Some definitions}\label{s601191}
\begin{itemize}
\item
Let $p$  be an odd prime.

Let $\zeta$  be a primitive $p$-th root of unity.

Let $K = \Q(\zeta)$  be the $p$-cyclotomic field.

Let $C^-$ be the relative class group of $K$.

Let $h$ be a prime, let us note $C(h)$ the $p$-Sylow subgroup of $C^-$.

Let $h^-(K)$ be the relative class number of $K$.

Let $K^+$ be the maximal totally real subfield of $K$.

Let $h(K^+)$ be the class group of $K^+$.

Let $G$ be the Galois group of the extension $K/\Q$.

Let $v$ be a primitive root $\modu p$.

Let $\sigma$ be the generator of $G$ given by $\sigma : \zeta \rightarrow\zeta^v$.

 For $n \in \Z$,
let us note $v_n$ for $v_n = v^n \modu p$  with $1 \leq v_n \leq p-1$.
\item
Let, for the indeterminate $X$, $P(X) \in \Z[X]$ given by $P(X) =
\sum_{i=0}^{p-2} v_{-i} X^i$. From Stickelberger relation, it can be
shown that $P(\sigma)$ annihilates $C^-$.

It can be also shown that for the indeterminate $X$,
\bd
P(X)(X-v) = p\times
Q(X)+v (X^{p-1}-1),
\ed
 where $Q(X) \in \Z[X]$   is given by $Q(X)=\sum_{i=1}^{p-2} \frac{v_{i-1}-v_{-i} v}{p}\times X^i$
and $Q(\sigma)$ annihilates the relative p-class
group $C_p^-$ of $K$ (the subgroup of exponent $p$ of $C^-$).
It can be shown that $\theta=\frac{P(\sigma)}{p}=\frac{Q(\sigma)}{\sigma-v}$ where  $\theta$ is  the Stickelberger element.
\end{itemize}
%
%
\section{Numerical observations}
\begin{itemize}
\item
I have observed, with a little MAPLE program that:

{\it For all the prime
numbers $p < 500$ and all
the odd prime numbers $h < p^2$ then the $h$-rank of the $h$-Sylow
subgroup $C(h)$ of $C^-$ is equal
to the degree of the polynomial $D(X)$ given by the following relations:
\bn
\item
If h not = p  then
\be\label{e702181}
D(X) = Gcd (P(X), X^{(p-1)/2} + 1)\modu h.
\ee
\item
if $h = p$ then
\be\label{e702182}
D(X) = Gcd (Q(X), X^{(p-1)/2} + 1) \modu h.
\ee
where
\be\label{e702183}
gcd(h, h^-(K)) = 1 \Leftrightarrow  degree(DX)) = 0.
\ee
\en}
\item
These observations are derived from a comparison of my results in Queme \cite{que} with R.
Schoof \cite{sch}  where the
structure of minus-class group for $p < 500$ is displayed.
\item
I think that,   for $h \not =p$  the  theorem 6.21 p. 103 in
Washington   \cite{was}
 and  for $h=p$ the relation, with Washington notations,  $[R_p^-:I_p]=p$-part of $h^-(K)$   in the end of the proof of theorem 6.21
p. 106  bring  a proof of  the following assertion :

\be\label{e702184}  degree
(D(X))   \not  = 0  \ \Leftrightarrow\    h\ |\ h^-(K),
\ee
 for  all
the odd prime numbers $h$,
assertion which generalizes my  numerical
observations in relation (\ref{e702183}).
Observe that this result uses only elementary arithmetical computations, hence is particularly suitable for the computer calculations
of the prime divisors  of the relative class number of the cyclotomic field $K$.
\item
Let $D(X)=\prod_k D_k(X)^{n_k}$ be the factorization of the polynomial $D(X)$ in ${\bf F}_h[X]$.  If $h=p$ then clearly $n_k=1$ for all $k$.
If $h\not = p$ there exist some primes $p$ with $n_k\not= 1$.
As an example for  $p=331, h=3$, then   $D(X)=(X+1)^2\times (X^4+2 X^3+X^2+2X+1)$ with the primitive root $v=3,
\sigma : \zeta\rightarrow \zeta^3$, and the rank of
$C(3)$ is $6$.
\end{itemize}
%
%
\section{A question  for the case $h\not= p $}
Observe that, unfortunately,  I am not  able to compute the structure of the class group $C^-$
to suppress the doubt  for  large primes  $p$ with my limited  {\it student} MAPLE software.
\paragraph{Question:}
Does there exists known counter-examples to my observations in relation (\ref{e702181}) p. \pageref{e702181}
when $p > 500$ and $h$ is
an odd prime with $h  \not = p$?
\paragraph{Comment:}
In regards to   the case $p = 4027$
and $h = 3$ candidate for a possible counter-example given
in Washington  p. 106, does   $R^-/I^-$ not isomorphic to $C^-$ implies
 $degree (D(X)) \not = 3$-rank of $C(3)$, which should imply that $(p=4027, h=3)$  is  a counter-example to my assertion
 (\ref{e702181})  p. \pageref{e702181}?

The case ($p=4027, h=3)$ described in Washington    is  not completely conclusive for our  question:
for instance we know from this description
that $R^-/I^- \simeq C^-$
is false  when the
 class group $C_2$  of the quadratic imaginary field  $\Q(\sqrt{-p})$ is not   cyclic :
 \bn
 \item
 Does it implies, looking at the full relative class group $C^-$, that
 the \underline{rank} of the $h$-Sylow subgroups of $R^-/I^-$ and of
 $C^-$ are different. Note that  $4026= 2.3.11.61$ and the difference of rank at quadratic field level  could perhaps disappear
 at the level of one of the intermediate fields between $\Q(\sqrt{-p})$ and $K$ or at the level of the field $K$ itself.
 \item
 Does it implies,  in  our elementary formulation,   that the rank  of $C(h)\not= degree(D(X))$?
 \en
 In the other hand, in answer to my  question, R. Schoof guess that the kind of behavior observed for quadratic fields
 also occurs for minus class groups of abelian fields of higher degree,
 more or less as the Cohen-Lenstra heuristics  predict and thus that counter-examples are likely.

%
\section{Context: }
The motivations   of this note  appear  more precisely in my
preprint entitled {\it On prime factors of class group of cyclotomic
fields} section 5 page 18. The URL of this  preprint is at :
\begin{verbatim}
http://arxiv.org/PS_cache/math/pdf/0609/0609723.pdf
\end{verbatim}
%

%
Roland Qu\^eme

13 avenue du ch\^ateau d'eau

31490 Brax

France

mailto: roland.queme@wanadoo.fr

home page: http://roland.queme.free.fr/

tel france 0561067020
\end{document}